\documentclass{article}
\usepackage{amssymb}
\usepackage[cp1251]{inputenc}
\usepackage[english]{babel}
\usepackage{amsmath}

\def\mapr#1{\smash{\mathop{\buildrel{#1}\over\longrightarrow}}}

\newlength{\myunit}
\newcounter{mycount}
\def\R{{\bf R}}

\def\C{{\bf C}}
\def\N{{\bf N}}

\def\K{{\bf K}}

\def\Fa{{\cal F}}

\def\Ga{{\cal G}}
\def\Ka{{\cal K}}
\def\Aa{{\cal A}}
\def\Ba{{\cal B}}
\def\Ea{{\cal E}}
\def\Da{{\cal D}}
\def\Ua{{\cal U}}
\def\Pa{{\cal P}}

\def\inxx#1{}
\def\id{{\bf id}}

\newtheorem{thm}{Theorem}
\newtheorem{lm}{Lemma}
\newtheorem{df}{Definition}

\newtheorem{cor}{Corollary}
\newtheorem{rem}{Remark}

\newcommand\qed{\hfill\vrule width2mm height2mm depth2mm}
\newcommand\proof{{\bf Proof. }\nobreak\noindent}

\newcommand\End{{\hbox{\bf End}\;}}

\newcommand\myker{{\hbox{\bf Ker}\;}}

\newcommand\im{{\hbox{\bf Im}\;}}

\def\mathrm#1{\hbox{\rm #1 }}

\title{On Compact and Fredholm Operators over C*-algebras and a New Topology
in the Space of Compact Operators}
\author{Anwar A. Irmatov and Alexandr S. Mishchenko}
\date{27th April 2005}
\begin{document}
\maketitle

\begin{abstract}
It is shown that the class of Fredholm operators
over an arbitrary unital $C^{*}$--algebra,
which may not admit adjoint ones, can be extended
in such a way that this class of
compact operators, used in the definition of
the class of Fredholm operators, contains
compact operators both with and without existence
of adjoint ones. The main property
of this new class is that a Fredholm operator which may
not admit an adjoint one has a
decomposition into a direct sum of an isomorphism and
a finitely generated operator.

In the space of compact operators in the Hilbert space
a new IM-topology is defined.
In the case when the $C^{*}$--algebra is a commutative
algebra of continuous functions on
a compact space the IM-topology fully describe
the set of compact operators over the
$C^{*}$--algebra without assumption of existence
bounded adjoint operators over the algebra.
\end{abstract}

\section{Introduction}
In the paper \cite{AtSe-03e} M.~Atiyah and G.~Segal have considered
families of Fredholm operators parametrized by points of a compact
space $K$ which are continuous in a topology weaker than the uniform
topology, i.e. the norm topology in the space of bounded operators
$B(H)$ in a Banach space $H$.

Therefore, it is interesting to ascertain whether the conditions,
characterized  families of Fredholm operators, from the paper
\cite{AtSe-03e} precisely describe the families of Fredholm operators
which forms a Fredholm operator over the $C^*$--algebra $\Aa=C(K)$ of
all continuous functions on $K$.

It is not supposed by the authors of the paper \cite{AtSe-03e}
that an operator over algebra $\Aa$ admits the adjoint one or in their
terms continuity of the adjoint family.

The aim of this paper is a fully clarification in the question of
description of the class of Fredholm operators which in general
case do not admit the adjoint operator.
For the first time, operators which play the role of Fredholm
operators and may not have the adjoint ones were considered in
the paper \cite{MiFo47}. Since the main class of operators considered
in the paper \cite{MiFo47} is the class of pseudodifferential operators,
for any element of which the adjoint operator automatically is a bounded
one, then existence of the adjoint operator was not the actual question
for the main goals of this paper.
However, in the paper \cite{AtSe-03e} authors have considered operators,
which may not have the adjoint one, in the form of families of operators
continuous in the compact-open topology the adjoint families of which, in
general case, may not be a continuous one.
In the present paper we show that the class of Fredholm operators over
arbitrary $C^*$--algebra, which may not admit the adjoint ones, can be
extended in a such way that the class of compact operators used in the
definition of the class of Fredholm operators contains compact operators
both with and without existence the adjoint ones.

In the case when the $C^*$--algebra is a commutative
algebra of continuous functions on a compact space
appropriate topologies in the classic spaces
of Fredholm and compact operators in the Hilbert space
are constructed
which fully describe the sets of Fredholm and compact
operators over the $C^*$--algebra without assumption of
existence bounded adjoint operators over the algebra.
A comparison with the class of operators
considered in the paper  \cite{AtSe-03e}  is given and it is shown
that the class of operators from the present paper strictly
includes the class of operators from the paper \cite{AtSe-03e}.

\section{A Notion of Compact Operator over $C^{*}$--algebra}

Let $A$ be a unital $C^{*}$--algebra. We shall consider so
called Hilbert $C^{*}$--modules
over the algebra $A$. The simplest Hilbert modules are the
free finitely generated $A$--modules
$$A^{n}=\underbrace{A\oplus A\oplus\cdots \oplus A}_{n
\mbox{ {\footnotesize   times}}}$$
and the $A$-module $l_{2}(A)=A^{\omega}$. All such modules have a
convenient description.
Any element $x$ of a module $A^{\alpha}$, $\alpha \in [1..\omega]$
is a sequence
$x=\{x_{1},x_{2},\dots,x_{k},\dots\}$, $x_{k}\in A$,
$1\leq k<1+\alpha$, such that the sum
\begin{equation}\label{eq3}
  \langle x,x\rangle=\sum\limits_{k=1}^{\alpha}x_{k}x_{k}^{*}\in A
\end{equation}
converges in the algebra $A$. It is clear that if $\alpha<\omega$
then the sum (\ref{eq3})
automatically converges.
The elements $e^{k}\in A^{\alpha}$, $e^{k}_{j}=\delta^{k}_{j}$ form a free
basis in the module $A^{\alpha}$ both for finite $\alpha$ and for infinite
$\alpha$, in the sense that any element $x\in A^{\alpha}$ can be represented
as a converged sum

\begin{equation}\label{eq4}
  x=\sum\limits_{k=1}^{\alpha}x_{k}e^{k}.
\end{equation}

 In general, a Hilbert $C^{*}$--module $M$
is a Banach space.
We say that $C^{*}$--module $M$ is a finitely generated $C^{*}$--module
if $M$ is a finitely generated $C^{*}$--module in the algebraic sense.
In other words,
there exists a free $C^{*}$--module $A^{n}$, $n<\omega$, and an algebraic
epimorphism
\begin{equation}\label{eq5}
f:A^{n}\mapr{}M\mapr{}0.
\end{equation}
It is easily verified that the epimorphism $f$ is a bounded map. Indeed,
if $x\in A^{n}$,
$x=\{x_{k}\}$ then
\begin{equation}\label{eq6}
  \begin{array}{c}
    \|f(x)\|^{2}=\|\langle f(x),f(x)\rangle\|=
    \left\|\left\langle \left(\sum\limits_{k} x_{k}f(e^{k})\right),
    \left(\sum\limits_{j} x_{j}f(e^{j})\right)\right\rangle\right\|=\\
     \\
    =\| \sum\limits_{k,j}\left(x_{k}\langle f(e^{k}),f(e^{j})\rangle x^{*}_{j}\right) \|\leq
     \| \sum\limits_{k,j}\|\langle f(e^{k}),f(e^{j})\rangle\|\cdot\left(x_{k} x^{*}_{j}\right) \|\leq\\ \\
\leq \sum\limits_{k,j}\|\langle f(e^{k}),f(e^{j})\rangle\|\cdot\left\|x_{k} x^{*}_{j}\right\|\leq
\sum\limits_{k,j}\|\langle f(e^{k}),f(e^{j})\rangle\|\cdot\|x_{k}\|\cdot\| x^{*}_{j}\|\leq
\\ \\
\leq\sum\limits_{k,j}\|\langle f(e^{k}),f(e^{j})\rangle\|\cdot\|x\|^{2}\leq
n^{2}C\|x\|^{2},
  \end{array}
\end{equation}
where
\begin{equation}\label{eq7}
  C=\mathop{\max}\limits_{k,j}\|\langle f(e^{k}),f(e^{j})\rangle\|.
\end{equation}

A Hilbert $C^{*}$--module is called a projective finitely generated
$C^{*}$--module if it is isomorphic to a direct summand of a finite
free $C^{*}$--module $L_{n}(A)=A^{n}$.

\begin{thm}\label{thm0} {\rm \cite{Mi}(Theorem 1.1, p. 69.) }
Let $M$ --- be a finitely generated Hilbert $A$--module.
Then $M$ is a projective $A$--module,
i.e. $M$ is isomorphic to a direct summand of a finite
free $A$--module $L_{n}(A)$.
\end{thm}

So, we can give the following definition

\begin{df}\label{df1}
Let $\End(l_{2}(A))$ be a Banach algebra of all bounded $A$--operators
of a Hilbert $A$--module $l_{2}(A)$. An $A$--operator
$K:l_{2}(A)\mapr{}l_{2}(A)$ is called a finitely generated
$A$--operator if it can be represented as a composition of
bounded $A$--operators $f_1$ and $f_2$:
$$K:l_{2}(A)\mapr{f_{1}}M\mapr{f_{2}}l_{2}(A),$$
where $M$ --- is a finitely generated Hilbert $C^{*}$--module.
The set $\Fa\Ga(A)\subset\End(l_{2}(A))$ of all finitely
generated $A$--operators forms a two side ideal. By definition,
an $A$--operator $K$ is called a compact if it belongs to the closure
$\Ka(l_2(A))=\overline{\Fa\Ga(A)}\subset\End(l_{2}(A))$, which also
forms two side ideal.
\end{df}

In general, the set $\Fa\Ga(A)\subset\End(l_{2}(A))$ is not closed subset.
For example, in classical case, when $\Aa=\C$, the set $\Fa\Ga(A)$
consists of all finite dimensional operators, while not all compact
operators are finite dimensional.

\begin{lm}\label{lm1}
The ideal $\Ka(l_2(A))$ is a proper ideal.
\end{lm}
\proof It is sufficient to prove that the identity
operator $\id\in\End(l_{2}(A))$ does not belong to $\Ka(l_2(A))$.
Or, to prove that the distance (in the sense of the operator norm)
between this operator and the set $\Fa\Ga(A)$ is a positive number.
In other words, it is sufficient to prove that any finitely generated
$A$--operator is not invertible. Indeed, if a finitely generated $A$--operator
$K:l_{2}(A)\mapr{f_{1}}M\mapr{f_{2}}l_{2}(A)$ is an invertible
$A$--operator then that means that the $A$--operator  $f_{2}$ is an
epimorphism. Since $C^{*}$--module $M$ is a finitely generated
$C^{*}$--module then there exists an epimorphism
$p:L_{n}(A)\mapr{}M$. Then the $A$--operator
$f_{2}\circ p:L_{n}(A)\mapr{}l_{2}(A)$ is an epimorphism. But this is
impossible.\qed

Let $l_{2}(A)= \left(L_{n}(A)\right)^{\perp} \oplus L_{n}(A)$ be an
orthogonal decomposition which is given by a pair of projectors
\begin{equation}\label{eq14a}
p_{n},q_{n}:l_{2}(A)\mapr{}l_{2}(A),
p_{n}+q_{n}=\id, \im p_{n} = L_{n}(A).
\end{equation}
Any $A$--operator
$f:l_{2}(A)\mapr{}l_{2}(A)$ forms a matrix composed from the bounded operators
\begin{equation}\label{eq14}
  f=\left(
  \begin{array}{cc}
    q_{n}fq_{n} & q_{n}fp_{n} \\
    p_{n}fq_{n} & p_{n}fp_{n} \
  \end{array}
  \right):\left(L_{n}(A)\right)^{\perp} \oplus L_{n}(A) \mapr{}
    \left(L_{n}(A)\right)^{\perp} \oplus L_{n}(A).
\end{equation}

\begin{thm}\label{thmCom} A bounded $A$--operator
$K:l_{2}(A)\mapr{}l_{2}(A)$ is a compact $A$-operator iff for
any $\varepsilon >0$ there exists a number $N$ such that for any $m>N$ we have
\begin{equation}\label{eqCom}
  \|q_{m}K\|\leq \varepsilon.
\end{equation}
\end{thm}
\proof Let us assume that the property (\ref {eqCom}) holds. Let $K_m=p_mK$. Since
\begin{equation}\label{eq16}
K_m:l_{2}(A)\mapr{f_{1}=p_mK}L_m(A)\mapr{f_{2}=i}l_{2}(A)
\end{equation}
then the operator $K_m$ is a finitely generated $A$--operator, i.e. $K_m\in \Fa\Ga(A)$. Since for any $\varepsilon>0$ there exists a natural number $N$ such that for any $m>N$
$$\|K-K_m\|=\|K-p_mK\|=\|q_mK\|\leq \varepsilon,$$
then $K\in\overline{\Fa\Ga(A)}$, i.e. the operator $K$ is a compact $A$-operator.

Inverse, Let $K$ be a compact $A$--operator. It follows from the definition \ref{df1} that there exists a finitely generated $A$--operator $K'\in \Fa\Ga(A)$ such that
\begin{equation}\label{eq15a}
  \|K-K'\|\leq \frac{\varepsilon}{2}.
\end{equation}

The finitely generated $A$--operator $K'$ can be represented as a composition
\begin{equation}\label{eq16a}
K':l_{2}(A)\mapr{f_{1}}M\mapr{f_{2}}l_{2}(A),
\end{equation}
in which, without loss of generality, we can assume that $M=L_{n}(A)$
with the basis
$e_{1},e_{2},\dots, e_{n}$. In the other words, the operator $f_{1}$ can
be described as linear combination of bounded functionals
\begin{equation}\label{eq17}
  f_{1}(x)=\sum\limits_{j=1}^{n}e_{j}\varphi^{j}(x), \quad \|\varphi_{j}\|\leq C.
\end{equation}
Correspondingly, the operator $f_{2}$ is given by a set of vectors
$y_{j}=f_{2}(e_{j})\in l_{2}(A)$. Thus, the operator $K'$ can be
represented by the formula
\begin{equation}\label{eq18}
  K'(x)=\sum\limits_{j=1}^{n}y_{j}\varphi^{j}(x).
\end{equation}
Under the formula (\ref{eq15a}) the operator $K'$ has the following matrix form:

\begin{equation}\label{eq19}
  K'=\left(
  \begin{array}{cc}
    q_{m}K'q_{m} & q_{m}K'p_{m} \\
    p_{m}K'q_{m} & p_{m}K'p_{m} \
  \end{array}
  \right):\left(L_{m}(A)\right)^{\perp} \oplus L_{m}(A) \mapr{}
    \left(L_{m}(A)\right)^{\perp} \oplus L_{m}(A).
\end{equation}

We have:
\begin{equation}\label{eq20}
q_{m}K'(x)= \sum\limits_{j=1}^{n}q_{m}(y_{j})\varphi^{j}\left(x\right).
\end{equation}
Then
\begin{equation}\label{eq20a}
\begin{array}{c}
 \|q_{m}K'(x)\|\leq \|\sum\limits_{j=1}^{n}q_{m}(y_{j})
 \varphi^{j}\left(x\right)\|\leq \\
 \leq\sum\limits_{j=1}^{n}\|q_{m}(y_{j})\|\cdot
 \|\varphi^{j}\|\cdot\|x\|.
\end{array}
\end{equation}
Since the number of vectors $y_{j}$ is finite then there exists
a number $N$ such that for any $m>N$
$\|q_{m}(y_{j})\|\leq \frac{\varepsilon}{2nC}$.
Then for any $m>N$ we have
\begin{equation}\label{eq21}
  \|q_{m}K'(x)\|\leq \frac{\varepsilon}{2}\|x\|,
\end{equation}
i.e. $\|q_{m}K'\|\leq \frac{\varepsilon}{2}$. Taking in account the inequality
(\ref{eq15a}) we obtain the desired inequality
\begin{equation}\label{eq23}
\|q_{m}K\|\leq \varepsilon.
\end{equation}
\qed
\vskip 0.2 cm

\begin{cor}\label{cor0}
Let $K:l_{2}(A)\mapr{}l_{2}(A)$ be a compact $A$-operator. Then for
any $\varepsilon >0$ there exists a number $N$ such that for any $m>N$ we have
\begin{equation}\label{eqCom1}
  \|q_{m}Kq_m\|\leq \varepsilon.
\end{equation}
\end{cor}
\proof We are interested in the operator $q_{m}K'q_{m}$ from the
formula (\ref{eq19}).

We have:
\begin{equation}\label{eq20}
q_{m}K'q_m(x)= \sum\limits_{j=1}^{n}q_{m}(y_{j})\varphi^{j}\left(q_m(x)\right).
\end{equation}
Then
\begin{equation}\label{eq20}
\begin{array}{c}
 \|q_{m}K'q_m(x)\|\leq \|\sum\limits_{j=1}^{n}q_{m}(y_{j})
 \varphi^{j}\left(q_m(x)\right)\|\leq \\
 \leq\sum\limits_{j=1}^{n}\|q_{m}(y_{j})\|\cdot
 \|\varphi^{j}\|\cdot\|q_m\|\cdot\|x\|.
\end{array}
\end{equation}
Since the number of vectors $y_{j}$ is finite then there exists a
number $N$ such that for any $m>N$
$\|q_{m}(y_{j})\|\leq \frac{\varepsilon}{2nC}$.
Then for any $m>N$ we have
\begin{equation}\label{eq21}
  \|q_{m}K'q_m(x)\|\leq \frac{\varepsilon}{2}\|x\|,
\end{equation}
i.e. $\|q_{m}K'q_m\|\leq \frac{\varepsilon}{2}$. Taking in account the inequality
(\ref{eq15a}) we obtain the desired inequality
\begin{equation}\label{eq23}
\|q_{m}Kq_m\|\leq \varepsilon.
\end{equation}
\qed
\vskip 0.2 cm

\section{Fredholm Operators over C*-algebra}

\begin{df}\label{df2}
A bounded $A$--operator $F:l_{2}(A)\mapr{}l_{2}(A)$ is called a
Fredholm $A$--operator if there exists
a bounded $A$--operator $G:l_{2}(A)\mapr{}l_{2}(A)$ such that
\begin{equation}\label{eq8}
  \id-FG\in\Ka(l_2(A)) , \quad \id-GF\in\Ka(l_2(A)).
\end{equation}
\end{df}

\begin{df}\label{df3}
We say that a bounded $A$--operator $F:l_{2}'(A)\mapr{}l_{2}''(A)$
admits an inner (Noether) decomposition if there is a decomposition
of the preimage and the image
\begin{equation}\label{eq9}
  l_{2}'(A)=M_{1}\oplus N_{1}, \quad l_{2}''(A)=M_{2}\oplus N_{2},
\end{equation}
where $C^{*}$--modules $N_{1}$ and $N_{2}$ are finitely generated
Hilbert $C^{*}$--modules, and if $F$ has the following matrix form
\begin{equation}\label{eq10}
F=\left(
  \begin{array}{cc}
    F_{1} & F_{2} \\
    0 & F_{4} \
  \end{array}
  \right):M_{1}\oplus N_{1}\mapr{}M_{2}\oplus N_{2},
\end{equation}
where $F_1:M_1\mapr{} M_2$ is an isomorphism.
\end{df}

\begin{df}\label{df31}
We put by definition $index \ F = [N_2]-[N_1]\in K(A)$.
\end{df}

\begin{df}\label{df4}
We say that a bounded $A$--operator $F:l_{2}'(A)\mapr{}l_{2}''(A)$
admits an external (Noether) decomposition if there exist finitely
generated $C^{*}$--modules $X_{1}$ and $X_{2}$ and bounded $A$--operators $E_{2}$, $E_{3}$ such that the matrix operator
\begin{equation}\label{eq11}
F_{0}=\left(
  \begin{array}{cc}
    F & E_{2} \\
    E_{3} & 0 \
  \end{array}
  \right):l_{2}'(A)\oplus X_{1}\mapr{}l_{2}''(A)\oplus X_{2}
\end{equation}
is an invertible operator.
\end{df}

\begin{df}\label{df41}
We put by definition $index \ F =[X_1]-[X_2]\in K(A)$.
\end{df}

\begin{thm}\label{thm1}
A bounded $A$--operator $F:l_{2}'(A)\mapr{}l_{2}''(A)$ admits an
external (Noether) decomposition iff
it admits an inner (Noether) decomposition.
\end{thm}
\proof If we have an inner (Noether) decomposition (\ref{eq10}) then
we can construct an external decomposition by an $A$--operator $F_{0}$
which has the following matrix form
\begin{equation}\label{eq12}
F_{0}=
\left(
  \begin{array}{ccc}
    F_{1} & F_{2} & 0 \\
    0 & F_{4} & \id \\
    0 & \id & 0 \
  \end{array}
  \right):M_{1}\oplus N_{1}\oplus N_{2}\mapr{}M_{2}\oplus N_{2}\oplus N_{1}.
\end{equation}
It is obvious that the operator $F_{0}$ is an invertible $A$--operator.

Now, let an external decomposition (\ref{eq11}) is given. Then the operator
$E_{3}:l_{2}'(A)\mapr{}X_{2}$ is an epimorphism.
Since the module $X_{2}$ is a projective $C^{*}$--module then there
exists a decomposition
\begin{equation}\label{eq12a}
  l_{2}'(A)=M_{1}\oplus N_{1}, \quad M_{1}=\myker E_{3},
  \quad E_{3}'=(E_{3})_{|N_{1}}:N_{1}\approx X_{2}.
\end{equation}

Analogously, let the inverted operator $G_{0}=F_{0}^{-1}$ has the
following matrix form
\begin{equation}\label{eq12b}
G_{0}=\left(
  \begin{array}{cc}
    G & G_{2} \\
    G_{3} & G_{4} \
  \end{array}
  \right):l_{2}''(A)\oplus X_{2}\mapr{}l_{2}'(A)\oplus X_{1}.
\end{equation}
The condition $G_{0}=F_{0}^{-1}$ can be rewrited as
$F_{0}G_{0}=\id_{(l_{2}''(A)\oplus X_{2})}$,
$G_{0}F_{0}=\id_{(l_{2}'(A)\oplus X_{1})}$, which have the following
matrix forms
\begin{equation}\label{eq12ca}
F_{0}G_{0}=
  \left(
  \begin{array}{cc}
    F & E_{2} \\
    E_{3} & 0 \
  \end{array}
  \right)
  \left(
  \begin{array}{cc}
    G & G_{2} \\
    G_{3} & G_{4} \
  \end{array}
  \right)=
  \left(
  \begin{array}{cc}
    \id & 0 \\
    0 & \id \
  \end{array}
  \right),
\end{equation}

\begin{equation}\label{eq12cb}
G_{0}F_{0}=
  \left(
  \begin{array}{cc}
    G & G_{2} \\
    G_{3} & G_{4} \
  \end{array}
  \right)
  \left(
  \begin{array}{cc}
    F & E_{2} \\
    E_{3} & 0 \
  \end{array}
  \right)=
  \left(
  \begin{array}{cc}
    \id & 0 \\
    0 & \id \
  \end{array}
  \right).
\end{equation}

The conditions (\ref{eq12ca}), (\ref{eq12cb}) can be rewrited as
\begin{equation}\label{eq12c}
  \begin{array}{rcl}
  \id=FG+E_{2}G_{3}&:& l_{2}''(A)\mapr{}l_{2}''(A);\\
  0=FG_{2}+E_{2}G_{4}&:& X_{2}\mapr{}l_{2}''(A);\\
  0=E_{3}G&:& l_{2}''(A)\mapr{}X_{2};\\
  \id=E_{3}G_{2}&:& X_{2}\mapr{}X_{2};\\
    \id = GF+G_{2}E_{3}&:&l_{2}'(A)\mapr{}l_{2}'(A);\\
    0=GE_{2}&:&X_{1}\mapr{}l_{2}'(A);\\
    0=G_{3}F+G_{4}E_{3}&:&l_{2}'(A)\mapr{}X_{1};\\
    \id=G_{3}E_{2}&:&X_{1}\mapr{}X_{1}. \
  \end{array}
\end{equation}
In particular, the operator $G_{3}:l_{2}''(A)\mapr{}X_{1}$ is also an
epimorphism.
Hence, there exists a decomposition
\begin{equation}\label{eq12a}
  l_{2}''(A)=M_{2}\oplus N_{2}, \quad M_{2}=\myker G_{3},\quad
  G_{3}'=(G_{3})_{|N_{2}}:N_{2}\approx X_{1}.
\end{equation}

Then the operator $F_{0}$ has the following matrix form
\begin{equation}\label{eq13}
  F_{0}=\left(
  \begin{array}{ccc}
    F_{1} & F_{2} & * \\
    0 & F_{4} & * \\
    0 & E_{3}' & 0 \
  \end{array}
  \right):M_{1}\oplus N_{1}\oplus X_{1}\mapr{}M_{2}\oplus N_{2}\oplus X_{2}.
\end{equation}

Indeed, if $x\in M_1$ then $E_{3}(x)=0$. Hence,
$G_{3}F(x)=0$, i.e. $F(x)\in \myker G_{3}=M_{2}$, and
$F_{0}(x)\in M_{2}$.
If $y\in M_{2}$ then $G_{3}(y)=0$, and $E_{3}G(y)=0$, i.e.
$G(y)\in M_{1}$. Moreover, if $x\in M_{1}$ then $x=GF(x)$, and for
$y\in M_{2}$ we have $y=FG(y)$. Hence, the operator $F_{1}$ is an
invertible $A$--operator.
Since the operators $E_{3}'$  and $G_{3}'$ are invertible $A$--operators then
the modules $N_{1}$ and $N_{2}$ are finitely generated Hilbert $C^*$--modules.
\qed

\begin{cor}\label{cor1}
The index constructed by inner or external decomposition does not depend
on the method of decomposition.
\end{cor}

\begin{thm}\label{thm2} Let $K:l_{2}(A)\mapr{}l_{2}(A)$ ---
be a compact operator in the sense of definition \ref{df1}.
Then the operator $\id + K$ admits an inner (Noether) decomposition.
\end{thm}
\proof Under the formula (\ref{eq14}) any operator $f:l_{2}(A)\mapr{}l_{2}(A)$
has the following matrix form:
\begin{equation}\label{eq14b}
  f=\left(
  \begin{array}{cc}
    q_{n}fq_{n} & q_{n}fp_{n} \\
    p_{n}fq_{n} & p_{n}fp_{n} \
  \end{array}
  \right):\left(L_{n}(A)\right)^{\perp} \oplus L_{n}(A) \mapr{}
    \left(L_{n}(A)\right)^{\perp} \oplus L_{n}(A).
\end{equation}

Due to the corollary \ref{cor0} we can find a natural number $N$ such
that for any $m>N$
\begin{equation}\label{eq15b}
  \|q_{m}Kq_{m}\|<1.
\end{equation}
The operator $F=\id+K$ can be represented in the following matrix form
\begin{equation}\label{eq24}
F=\left(
  \begin{array}{cc}
    F_{1} & F_{2} \\
    F_{3} & F_{4} \
  \end{array}
  \right):
  \left(L_{m}(A)\right)^{\perp} \oplus L_{m}(A) \mapr{}
    \left(L_{m}(A)\right)^{\perp} \oplus L_{m}(A),
\end{equation}
where the operator $F_{1}$ has the form $F_{1}=\id+q_{m}Kq_{m}$, and hence,
is an invertible $A$-operator.  The invertibility of the operator $F_{1}$
allows to represent the matrix (\ref{eq24}) in the following form
\begin{equation}\label{eq25}
\begin{array}{c}
  \left(
  \begin{array}{cc}
    F_{1} & F_{2} \\
    F_{3} & F_{4} \
  \end{array}
  \right)= \\ \\ =
\left(
  \begin{array}{cc}
    \id & 0 \\
    F_3F_1^{-1}0 & \id \
  \end{array}
  \right)\cdot
\left(
  \begin{array}{cc}
    F_{1}& 0 \\
     0& F_{4}-F_3F_1^{-1}F_2 \
  \end{array}
  \right)\cdot
\left(
  \begin{array}{cc}
    \id & F_1^{-1}F_2 \\
     0 & \id \
  \end{array}
  \right),
\end{array}
\end{equation}
This proves the theorem.
\qed

\begin{thm}\label{thm3}
Any Fredholm operator in the sense of definition \ref{df2} admits
both the inner and external (Noether) decomposition.
\end{thm}
\proof Let operators $F:l_{2}'(A)\mapr{}l_{2}''(A)$,
$G:l_{2}''(A)\mapr{}l_{2}'(A)$ are chosen such that
\begin{equation}\label{eq26}
  K'=\id-FG\in\Ka(l_2(A)) , \quad K''=\id-GF\in\Ka(l_2(A)).
\end{equation}
In accordance with the theorem \ref{thm2} there exist decompositions
\begin{equation}\label{eq27}
  \begin{array}{c}
    l_{2}'(A)=M_{1}\oplus N_{1}, \quad M_{1}=\im p_{1},
    \quad N_{1}=\im(\id-p_{1}),\\
    l_{2}'(A)=M_{2}\oplus N_{2}, \quad M_{2}=\im p_{2},
    \quad N_{2}=\im(\id-p_{2}),\
  \end{array}
\end{equation}
such that the modules $N_{1}$ and $N_{2}$ are finitely
generated $C^*$--modules, and the matrix of the operator
$\id-K''$ has a diagonal form
\begin{equation}\label{eq28}
  \id-K''=GF=\left(
  \begin{array}{cc}
    K_{1} & 0 \\
    0 & K_{2} \
  \end{array}
  \right):M_{1}\oplus N_{1}\mapr{F}l_{2}''(A)\mapr{G}M_{2}\oplus N_{2},
\end{equation}
where the operator $K_{1}$ is invertible. Let us consider the operator
\begin{equation}\label{eq29}
  P:l_{2}''(A)\mapr{}l_{2}''(A),\quad P(x)= FK_{1}^{-1}p_{2}G(x).
\end{equation}
We have
\begin{equation}\label{eq30}
\begin{array}{c}
   PP(x)=FK_{1}^{-1}p_{2}G\cdot FK_{1}^{-1}p_{2}G(x)=
FK_{1}^{-1}p_{2}K_{1}K_{1}^{-1}p_{2}G(x)=\\ \\ =
  FK_{1}^{-1}p_{2}p_{2}G(x)=FK_{1}^{-1}p_{2}G(x)=P(x),
\end{array}
\end{equation}
i.e. the operator $P$ is a projector. This means that the
module $l_{2}''(A)$ can be decomposed in direct sum
\begin{equation}\label{eq31}
  l_{2}''(A)=\im P\oplus \myker P=M_{3}\oplus N_{3},
\end{equation}
and in the decomposition (\ref{eq31}) the operator $F$ has
the following matrix form
\begin{equation}\label{eq32}
  F=\left(
  \begin{array}{cc}
    F_{1} & * \\
    0 & F_{4} \
  \end{array}
  \right):M_{1}\oplus N_{2}\mapr{}M_{3}\oplus N_{3},
\end{equation}
where the operator $F_{1}$ is an isomorphism.

Now it is necessary to prove that the module $N_{3}$ is a finitely
generated $C^*$--module. For, by the theorem \ref{thm2} the
operator $K'=\id -FG$ in the decompositions
\begin{equation}\label{eq34}
  \begin{array}{c}
    l_{2}''(A)=M_{4}\oplus N_{4},\\
    l_{2}''(A)=M_{5}\oplus N_{5} \
  \end{array}
\end{equation}
has the following matrix form
\begin{equation}\label{eq35}
\id-K'=FG:\left(
\begin{array}{cc}
  K_{3} & 0 \\
  0 & K_{4}
\end{array}
\right):M_{4}\oplus N_{4}\mapr{}M_{5}\oplus N_{5},
\end{equation}
where the operator $K_{3}$ is an isomorphism. In particular, the operator
\begin{equation}\label{eq36}
A=\left(
\begin{array}{ccc}
  K_{3} & 0 &0\\
  0 & K_{4}&\id
\end{array}
\right):M_{4}\oplus N_{4}\oplus N_{5}\mapr{}M_{5}\oplus N_{5}
\end{equation}
is an epimorphism. It is convenient to represent this operator
in the following matrix form
\begin{equation}\label{eq37}
A=\left(
\begin{array}{ccc}
  FG& , &a
\end{array}
\right):l_{2}''(A)\oplus N_{5}\mapr{}l_{2}''(A).
\end{equation}
We can represent the operator $A$ as composition
\begin{equation}\label{eq38}
A=
\left(
\begin{array}{ccc}
  F& , &a
\end{array}\right)
\left(
\begin{array}{cc}
  G &0 \\
  0&\id
\end{array}
\right):l_{2}''(A)\oplus N_{5}\mapr{}\l_{2}'(A)\oplus N_{5}\mapr{}l_{2}''(A).
\end{equation}
Hence, the operator
\begin{equation}\label{eq39}
B=
\left(
\begin{array}{ccc}
  F& , &a
\end{array}\right)
:l_{2}'(A)\oplus N_{5}\mapr{}l_{2}''(A)
\end{equation}
also is an epimorphism.

By (\ref{eq32}) the operator $B$ has the
following matrix form
\begin{equation}\label{eq40}
  B=\left(
  \begin{array}{ccc}
    F_{1} & *&a_{1} \\
    0 & F_{4}&a_{2} \
  \end{array}
  \right):M_{1}\oplus N_{2}\oplus N_{5}\mapr{}M_{3}\oplus N_{3}.
\end{equation}
Hence the operator
\begin{equation}\label{eq41}
  D=\left(
  \begin{array}{ccc}
  F_{4}&,&a_{2} \
  \end{array}
  \right):N_{2}\oplus N_{5}\mapr{}N_{3}.
\end{equation}
is an epimorphism. This means that the module $N_{3}$ is a
finitely generated $C^*$--module.
\qed
\vskip 0.2 cm

\begin{cor}\label{cor1a} Let $K:l_{2}(A)\mapr{}l_{2}(A)$ ---
be a compact operator in the sense of definition \ref{df1} and
$F:l_{2}(A)\mapr{}l_{2}(A)$ be a Fredholm $A$--operator in the
sense of definition \ref{df2}.
Then the operator $F + K$ is a Fredholm $A$--operator in the sense
of definition \ref{df2} and admits an inner (Noether) decomposition.
\end{cor}

\section{Fredholm and Compact Operators over
Commu\-tative $C^*$--algebras}

Let $\Aa=C(X)$ be an algebra of continuous functions on a compact
Hausdorff space $X$. We can identify the Hilbert $\Aa$--module
$l_2(\Aa)$ with the set $[X, H]$ of all continuous maps from the
space $X$ into the Hilbert space $H$. Denote the set of all maps
from the space $X$ into the space of bounded linear operators $\Ba(H)$,
continuous in the strong topology, by $[X, \Ba(H)^s]$. We denote
by $\Ea$ the map
\begin{equation}\label{eq42}
\Ea: [X, \Ba(H)^s] \to End_{\Aa}(l_2(\Aa))
\end{equation}
defined by the formula
\begin{equation}\label{eq43}
(\Ea(T)(\varphi))(x) = T(x)\varphi (x),
\end{equation}
where $T\in [X, \Ba(H)^s]$ is fixed, $x\in X$, and
$\varphi \in l_2(\Aa)=[X, H]$ are arbitrary,
and denote by $\Da$ the map
\begin{equation}\label{eq44}
\Da: End_{\Aa}(l_2(\Aa)) \to [X, \Ba(H)^s]
\end{equation}
defined be the formula
\begin{equation}\label{eq45}
(\Da(B)(x))(a)=(B\varphi )(x),
\end{equation}
where $B\in End_{\Aa}(l_2(\Aa))$ and the map $\varphi \in [X, H]$ is
chosen so that $\varphi (x) \equiv a$, $a\in H$.

It was shown in the paper \cite{Frank} that the definitions
(\ref{eq42})--(\ref{eq45}) are correct,
the maps $\Ea$
and $\Da$ are $\Aa$-isomorphisms and
\begin{equation}\label{eq46}
\Ea \Da = \id_{End_{\Aa}(l_2(\Aa))} \quad \Da \Ea = \id_{[X, \Ba(H)^s]}.
\end{equation}
Further, it was shown ibidem that each invertible $\Aa$--operator
$S\in GL(l_2(\Aa))$ under the map $\Da$ can be represented as a family
of invertible operators $S_x \in GL(H)^s$ continuous in the strong topology
such that $\sup_{x\in X}\|S_x^{-1}\|<\infty$, and, conversely, each family
of invertible operators $S_x \in GL(H)^s$ continuous in the strong topology
such that $\sup_{x\in X}\|S_x^{-1}\|<\infty$ is mapped by $\Ea$ into an
invertible $\Aa$--operator.

In the paper \cite{Irm} the F-topology was introduced in the space of
Fredholm operators $\Fa(H)$.

\begin{df}\label{F-top} {\rm \cite{Irm}} The following sets form a subbase
of the $F$-topology
$$
U_{\varepsilon,a_1,\ldots,a_n,A}=\left\{ B\in \Fa(H)\ |\
\Vert(B-A)a_i\Vert < \varepsilon \ \ \forall \  i=1,\ldots, n \right\},
$$
$$
U_{\varepsilon, V, A}= \left\{ B\in \Fa(H)\ |\ \exists R\in GL(H),
\ R(V)\subset V, \ such \  that \ \Vert RB-A\Vert < \varepsilon \right\}.
$$
Here $V$ denotes a finite dimensional subspace of the Hilbert space $H$
and  $a_1,\ldots,a_n \in H$.
\end{df}

Let $f:[0,1]\to \Fa(H)^F$ be any continuous map in the $F$--topology.
Then, $index\ f(x) = const$. On the other hand, there exists a
map $f:[0,1]\to \Fa(H)^s$ continuous in the strong topology such
that $index\ f(0)\ne index\ f(1)$ (see \cite{Irm}), so the $F$--topology
is strictly stronger than the strong topology in the space of Fredholm operators.

Let $F\in End_{\Aa}(\l_2(\Aa))$ be any Fredholm $\Aa$--operator.
Then for any $x\in X$ $(\Da(F))(x)\in \Fa(H)$ and the map
$\Da(F):X\to \Fa(H)^s$ is continuous in the strong topology.
It was shown in \cite{Irm} that the map
$$\Da(F):X\to \Fa(H)^F \stackrel{\id}{\longrightarrow} \Fa(H)^s
\subset \Ba(H)^s$$
is continuous in the $F$--topology and vice versa if a map
$f:X\to \Fa(H)^F$ is continuous in the $F$--topology then the
$\Aa$--operator $\Ea(f)$ is a Fredholm $\Aa$--operator. Thus, the map
\begin{equation}\label{eq47}
\Da|_{\Fa(\l_2(\Aa))}: \Fa(l_2(\Aa)) \to [X, \Fa(H)^F],
\end{equation}
where $\Fa(l_2(\Aa))$ is the space of Fredholm operators
over the algebra $\Aa$ and  $[X, \Fa(H)^F]$ is the set
of continuous maps  from the space $X$ into the space
of Fredholm operators $\Fa(H)^F$, with the $F$-topology,
is an isomorphism.

Denote by $\Ua(H)^s$ the space of unitary operators in $H$ with
the strong topology. Due to the formula
\begin{equation}\label{eq47a}
U(x)^{-1} -U(x_0)^{-1} = U(x)^{-1}\left(U(x_0) - U(x)\right) U(x_0)^{-1},
\end{equation}
we can assert that if a map $U:X\to \Ua(H)^s$ is continuous
in the strong topology then the map $U^{-1}:X\to \Ua(H)^s$ is
also continuous in the strong topology.

\begin{thm}\label{thm6}
Let $X$ be a compact Hausdorff space and maps
$U:X \to \Ua(H)^s$, $F:X\to \Fa(H)^F$ are continuous in the strong
topology and $F$-topology, respectively. Then the map
$UFU^{-1}:X\to \Fa(H)^s\subset \Ba(H)^s$, given by the
formula $UFU^{-1}(x)=U(x)F(x)U^{-1}(x)$, is continuous in the $F$-topology:
$$UFU^{-1}:X\to \Fa(H)^F \stackrel{\id}{\longrightarrow}\Fa(H)^s\subset
\Ba(H)^s.$$
\end{thm}

\proof Since the $\Aa$--operators $\Ea(U)$, $\Ea\left(U^{-1}\right)$
are unitary $\Aa$--operators, and $\Ea(F)$ is Fredholm $\Aa$--operator
then the operator $\Ea\left(UFU^{-1}\right)=\Ea(U)\Ea(F)\Ea\left(U^{-1}\right)$
is Fredholm $\Aa$--operator. Hence, by the isomorphism (\ref{eq47})
the map $UFU^{-1}=\Da\Ea\left(UFU^{-1}\right)$ is continuous in the $F$-topology.
\qed

Let us consider the set of compact $\Aa$--operators
$\Ka(l_2(\Aa))$. In the paper \cite{MiFo47} has been considered
the following class of compact operators $\Ka^*(l_2(\Aa))$. By
definition (see \cite{MiFo47}) an
$\Aa$--operator $K: l_2(A) \to l_2(A)$ belongs to the
set $\Ka^*(l_2(\Aa))$ iff
\begin{equation}\label{eq49}
\lim_{n\to \infty}\|Kq_n\| =0,
\end{equation}
where the operator $q_n$ is defined by the formula (\ref{eq14a}).
It was shown in (\cite{MaTr-01}, Prop. 2.2.1.) that the set $\Ka^*(l_2(\Aa))$
coincides with the closure of the set of linear combinations of elementary
operators $\theta_{x,y}(z):=x<y,z>$, where $x,y,z \in l_2(\Aa)$.
Hence, any $K\in \Ka^*(l_2(\Aa))$ automatically admits the adjoint operator.
On the other hand, our notion of compact
$\Aa$--operator does not demand existing of adjoint operator
unlike that was assumed in many papers on $KK$--theory and so we
shall distinguish the set $\Ka(l_2(\Aa))$
of all compact $\Aa$--operators
and the subset
$\Ka^*(l_2(\Aa))\subset \Ka(l_2(\Aa))$
of compact $\Aa$--operators which admit adjoint operator.

\begin{thm}\label{thm7} {\rm \cite{Irm}}
A compact $\Aa$--operator $K$ admits adjoint operator,
i.e. $K\in \Ka^*(l_2(\Aa))$,
iff the map
$$\Da(K): X \to \Ka(H)^u \stackrel{\id}{\longrightarrow} \Ka(H)^s
\subset \Ba(H)^s$$
is continuous in the uniform topology.
\end{thm}

The following example shows that there exists a self-adjoint family of
compact operators continuous
in the strong topology such that the corresponding
$\Aa$--operator does not belong to the set $\Ka(l_2(\Aa))$.

\vskip 0.2 cm
\noindent
{\bf Example.} Let $X=\{0\}\cup \bigcup_{i=1}^{\infty}\{\frac{ 1}{i}\}
\subset \R$.
We define the map $K:X\to \Ka(H)$ by the following formula
\begin{equation}\label{eq50}
K\left(\frac{1}{i}\right)(\xi)=-\xi_i,\qquad
K(0)=0,
\end{equation}
where $\xi=(\xi_1,\xi_2,\ldots)$ is an element of the standard Hilbert space.
Then $K\left(\frac{1}{i}\right)(\xi)\to 0$ as $i\to \infty$ for
every $\xi \in H$. But $\Ea(K) \notin \Ka(l_2(\Aa))$. Indeed, if we
suppose the contrary then by the theorem \ref{thm2} the operator
$\id +\Ea(K)$ is a Fredholm $\Aa$--operator. Due to the isomorphism
(\ref{eq47}) the map $\Da(\id +\Ea(K))=\id + K:X\to \Fa(H)^F$ is
continuous in the $F$-topology. But for any invertible operator $S$
we have
$$\left\|S\left(\id+ K\left(\frac{1}{i}\right)\right) - (\id + K(0))\right\| \geq \left\|\left(S\left(\id+ K\left(\frac{1}{i}\right)\right) - \id\right)(e_i)\right\|=1.$$
That means that the map $\id + K:X\to \Fa(H)^F$ is not continuous in
the $F$-topology.
\qed

Thus, the example poses the problem of finding a topology in the
space of compact operators $\Ka(H)$ such that any family continuous
in this topology forms a compact $\Aa$--operator, and vice versa,
any compact $\Aa$--operator maps by the map $\Da$ to a family of
compact operators continuous in the sought topology.

We define a new $IM$-topology in the space of compact operators
$\Ka(H)$ in the following way.
Let
$$
U_{\varepsilon,a_1,\ldots,a_n,K}=\left\{ B\in \Ka(H)\ |\
\Vert(B-K)a_i\Vert < \varepsilon \ \ \forall \  i=1,\ldots, n \right\},
$$
$$
U_{\varepsilon,n, S, K}= \{ B\in \Ka(H)\ |\  \exists R\in GL(H),
\ such \ that
$$
$$
\ \Vert R(S+Q_nB)-(S+Q_nK)\Vert < \varepsilon \},
$$
where $\varepsilon >0$, $S\in GL(H)$, and
$Q_n: H=\left(L_{n}\right)^{\perp} \oplus L_{n} \mapr{}
\left(L_{n}\right)^{\perp} \subset H$ is the orthogonal projection
along the subspace $L_n$ spanned by the first $n$ orthonormal
basis vectors $e_1, \ldots, e_n$.

\begin{df}\label{IM-top}
As a subbase of the $IM$-topology we take the following sets
$$U_{\varepsilon,a_1,\ldots,a_n,K} \quad and \quad U_{\varepsilon, S, K}
:=\bigcap_{n=0}^{\infty}U_{\varepsilon,n, S, K}.$$
\end{df}

\begin{rem}\label{rem1}
It follows from the definition of $IM$--topology that the identity map
$$\Ka(H)^{IM} \stackrel{\id}{\longrightarrow} \Ka(H)^s \subset \Ba(H)^s$$
from the space of compact operators with the $IM$--topology to the
same space with the strong topology is continuous. Since any sets
$U_{\varepsilon,a_1,\ldots,a_n,K} \ and \ U_{\varepsilon, S, K}$ contain
the ball $B(K,\varepsilon)=\{Z\in \Ka(H)| \ \|Z-K\|<\varepsilon\}$, then
the map
$$\Ka(H)^u \stackrel{\id}{\longrightarrow} \Ka(H)^{IM}$$
from the space of compact operators with the norm topology to the
same space with the $IM$--topology is continuous.
\end{rem}

\begin{thm}\label{thm8}
An $\Aa$--operator $K$ is compact operator, i.e. $K\in \Ka(l_2(\Aa))$,
iff the map
$$\Da(K): X \to \Ka(H)^{IM}\stackrel{\id}{\longrightarrow}
\Ka(H)^s \subset \Ba(H)^s$$
is continuous in the $IM$-topology.
\end{thm}
\proof Let $K\in \Ka(l_2(\Aa))$. We have to prove that the map $\Da(K)$
is continuous map from $X$ to $\Ka(H)$ with the $IM$--topology.
For, it is sufficient to show that for any $\varepsilon>0$, for
any $S\in GL(H)$, and for any $x_0\in X$ there exists
a neighbourhood $U_{x_0}\subset X$, $x_0\in U_{x_0}$ such that for
any $n\geq 0$
\begin{equation}\label{eq5010}
\Da(K)\left(U_{x_0}\right) \subset U_{\varepsilon, n, S, \Da(K)(x_0)}.
\end{equation}

Let $l_{2}(\Aa)= L_{n, \Aa}^{\perp} \oplus L_{n,\Aa}$ be an
orthogonal decomposition which is given by a pair of projectors
$p_n$, $q_n$, $p_n + q_n =\id$, $\im p_{n} = L_{n,\Aa}$,
$\im q_{n} = L_{n, \Aa}^{\perp}$, $\Da(q_n)(x)\equiv Q_n$.

Let $s:X\to GL(H)$ be a constant map, $s(x) \equiv S \in GL(H)$,
and $\hat S= \Ea(s)\in GL(\Aa)$, i.e. $\Da(\hat S)(x) \equiv S \in GL(H)$.

Let us choose $n\in \N$ such that for all $m>n$ the
$\Aa$--operator $G_m = \hat S + q_mK$ is invertible, i.e.
\begin{equation}\label{eq5011}
G_m = \hat S + q_mK \in GL(l_2(\Aa)).
\end{equation}
Then for any $x\in X$ we have
$$\id =\Da(G_m^{-1}G_m)(x)=\Da(G_m^{-1})(x)(S + Q_m\Da(K)(x)).$$
If we put $R:= \Da(G_m)(x_0)\Da(G_m^{-1})(x) \in GL(H)$ then for
any $x\in X$ and $m>n$ we have
$$
\|\Da(G_m)(x_0)\Da(G_m^{-1})(x)\left(S+Q_m\Da(K)(x)\right) -
\left(S+Q_m\Da(K)(x_0)\right)\|= 0.
$$
The last equality means that $\forall  x\in X$, $\forall m>n$,
and $\forall \varepsilon >0$
$$\Da(K)(x)\in U_{\varepsilon, m, S, \Da(K)(x_0)}.$$

Since $K$ is a compact $\Aa$--operator then by the corollary \ref{cor1a}
for any Fredholm $\Aa$--operator $F$ the operator $F+q_lK$ is
a Fredholm $\Aa$--operator. In particular, the operator $\hat S + q_lK$
is a Fredholm $\Aa$--operator for all $l\geq 0$. By the isomorphism
(\ref{eq47}), the maps $\Da(\hat S + q_lK):X\to \Fa(H)^F$,
$\Da(\hat S + q_lK)(x)= S + Q_l\Da(K)(x)$, $l=1,\ldots, n$, are
continuous in the $F$-topology.
Hence, for any fixed finite dimensional subspace $V_l\subset H$,
$l=1,\ldots, n$, there exists a neighbourhood $U_{x_0}^l\subset X$
such that for any $x\in U_{x_0}^l$
$$\Da(\hat S + q_lK)(x)\in U_{\varepsilon, V_l, S+Q_l\Da(K)(x_0)},$$
where $U_{\varepsilon, V_l, S+Q_l\Da(K)(x_0)}$ is an open set in $F$-topology.
This means that $\exists R\in GL(H), \ R(V_l)\subset V_l$, such that
$$\|R\left(S+Q_l\Da(K)(x)\right) - \left(S+Q_l\Da(K)(x_0)\right)\|
<\varepsilon,$$
i.e.
$$\Da(K)(x) \in U_{\varepsilon, l, S, \Da(K)(x_0)}.$$

Then, if we put $U_{x_0}:= \bigcap_{l=1}^n U_{x_0}^l$ we obtain
the necessary condition (\ref{eq5010}).
This proves that the map
$\Da(K):X\to \Ka(H)^{IM}\stackrel{\id}{\longrightarrow}
\Ka(H)^s \subset \Ba(H)^s$ is continuous in the $IM$-topology.

\vskip 0.2 cm

To prove the inverse assertion of the theorem let us suppose
the contrary. This means that there exists a continuous map
$C: X\to \Ka(H)^{IM}$ such that the operator $\Ea(C)$ is not
a compact $\Aa$-operator, i.e. $\Ea(C)\notin \Ka(l_2(\Aa))$.
Due to the theorem \ref{thmCom} there exist a number $c_1>0$
and an increasing sequence of natural numbers $\{n_i\}_{i\in \N}$ such that
\begin{equation}\label{eq50a}
\left\|q_{n_i}\Ea(C)\right\| > c_1.
\end{equation}
Since
\begin{equation}\label{eq50b}
\left\|q_{n_i}\Ea(C)\right\| = \sup_{x\in X}\left\|Q_{n_i}C(x)\right\|,
\end{equation}
then there exist an element $x_i\in X$ and a vector $v^i\in H$, $\|v^i\|=1$,
such that
\begin{equation}\label{eq50c}
\left\|Q_{n_i}C(x_i)(v^i)\right\| > c_1.
\end{equation}
Let
\begin{equation}\label{eq50cc}
c_0:=\frac{c_1}{2} \quad and \quad \varepsilon_0:=\min\left(c_0,
\frac{c_0^2}{\|\Ea(C)\|}\right).
\end{equation}
Since the map $C$ is continuous in the $IM$--topology which is stronger
than the strong topology then we can conclude that $\|\Ea(C)\| < \infty$
and $\varepsilon_0 >0$.

We assert that we can choose from the sequence $\{x_i\}_{i\in \N}$
a subsequence $\{y_i\}_{i\in \N}$, an increasing subsequence of natural
numbers $\{r_i\}_{i\in \N}$, and an orthonormal sequence of vectors
$w_i\in H$, $i\in \N$, such that
\begin{equation}\label{eq50d}
\left\|Q_{r_i}C(y_i)(w_i)\right\| > c_0
\end{equation}
and
\begin{equation}\label{eq50e}
\left\|Q_{r_i}C(y_j)(w_j)\right\| <\frac{\varepsilon_0}{2^{i+2}} \quad
for \  all \quad j<i.
\end{equation}
We shall prove our assertion by mathematical induction.
Let us put $y_1:=x_1$ and $w_1:= v^1$. Let us suppose that $k$
points $y_1,\ldots,y_k \in X$ and $k$ orthonormal vectors
$w_1,\ldots,w_k\in H$ have already been chosen such that the
conditions (\ref{eq50d}) and (\ref{eq50e}) hold. We consider $2k$
functions $\varphi_i:X\to H$, $\varphi_i(x)\equiv w_i$, and
$\Ea(C)\varphi_i$, $(\Ea(C)\varphi_i)(x)=C(x)w_i$, $i=1,\ldots,k$,
as elements of $\Aa$--module $l_2(\Aa)=[X,H]$. For any $\varepsilon >0$
there exists a natural number $N(\varepsilon)$ such that for any natural
number $l>N(\varepsilon)$ the inequality
\begin{equation}\label{eq50f}
\left\|q_l\Ea(C)\varphi_j\right\| < \varepsilon
\end{equation}
holds for $j=1,\ldots,k$.
Let $n_m>N\left(\frac{\varepsilon_0}{\sqrt{k}2^{k+3}}\right)$
and $H_1:= span<w_1,\ldots,w_k>$. Let $v^m=v^m_1 + v^m_2$,
$v^m_1\in H_1$, $v^m_2\in (H_1)^{\perp}$, be a representation of
the vector $v^m$ in accordance with the decomposition
$H=H_1\oplus (H_1)^{\perp}$.
We have

\begin{eqnarray}\label{eq50g}
c_1<\left\|\left(Q_{n_m}C(x_m)\right)(v^m)\right\| =
\left\|\left(Q_{n_m}C(x_m)\right)(v^m_1 + v^m_2)\right\|\leq && \nonumber \\
\frac{\sqrt{k}\varepsilon_0}{\sqrt{k}2^{k+3}}\left\|v^m_1\right\| +
\left\|\left(Q_{n_m}C(x_m)\right)(v^m_2)\right\|.&&
\end{eqnarray}
If we put $y_{k+1}:=x_m$, $w_{k+1}:=\frac{v^m_2}{\left\|v^m_2\right\|}$,
and $r_{k+1}:=n_m$, then by (\ref{eq50g}) we obtain the inequality
(\ref{eq50d}) for $i=k+1$:
\begin{equation}\label{eq50h}
\left\|\left(Q_{r_{k+1}}C(y_{k+1})\right)(w_{k+1})\right\|\geq
\left\|\left(Q_{n_m}C(x_m)\right)(v^m_2)\right\|>
c_1-\frac{\varepsilon_0}{2^{k+3}}>c_0.
\end{equation}
The inequality (\ref{eq50e}) for $i=k+1$ follows from the
inequality (\ref{eq50f}):
$$ \left\|\left(Q_{r_{k+1}}C(y_{j})\right)(w_j)\right\|
\leq \left\|q_{{n_m}}\Ea(C)\varphi_j\right\|\leq
\frac{\varepsilon_0}{\sqrt{k}2^{k+3}}, \quad j\leq k.$$

Let $a=(a_1,a_2,\ldots)\in H$, $\sum_{i=1}^{\infty}a_i\bar a_i=1$.
Let us estimate the norm of the element
$\sum_{i=s}^ta_i\cdot (Q_{r_i}C(y_i))(w_i)\in H$. We have
\begin{eqnarray}\label{eq50i}
\left\|\sum_{i=s}^ta_i\cdot (Q_{r_i}C(y_i))(w_i)\right\|^2 =
\left(\sum_{i=s}^ta_i\cdot (Q_{r_i}C(y_i))(w_i),
\sum_{i=s}^ta_i\cdot (Q_{r_i}C(y_i))(w_i)\right)= && \nonumber \\
\sum_{i=s}^t|a_i|^2\left(Q_{r_i}C(y_i)(w_i), Q_{r_i}C(y_i)(w_i)\right) &&
\nonumber \\
+ 2\sum_{s\leq i<j\leq t}{\bf Re}\left(a_i\cdot (Q_{r_i}C(y_i))(w_i), a_j
\cdot (Q_{r_j}C(y_j))(w_j)\right) \leq && \nonumber \\
\|\Ea(C)\|^2\sum_{i=s}^t|a_i|^2 + 2\sum_{s\leq i<j\leq t}{\bf Re}
\left(a_i\cdot (Q_{r_j}C(y_i))(w_i), a_j\cdot (Q_{r_j}C(y_j))(w_j)\right)
 \leq && \nonumber \\
\|\Ea(C)\|^2\sum_{i=s}^t|a_i|^2 + 2\|\Ea(C)\|\sum_{i=s}^{t-1}
\sum_{j=i+1}^t\frac{\varepsilon_0}{2^{j+2}} =
&& \nonumber \\
\|\Ea(C)\|^2\sum_{i=s}^t|a_i|^2 + 2\|\Ea(C)\|\varepsilon_0
\left(\frac{1}{2^{s+1}} -\frac{t-s+2}{2^{t+2}}\right),
&&
\end{eqnarray}
and
\begin{eqnarray}\label{eq50j}
\left\|\sum_{i=1}^ta_i\cdot (Q_{r_i}C(y_i))(w_i)\right\|^2 \geq && \nonumber \\
\left|c_0^2\sum_{i=1}^t|a_i|^2 -
2\sum_{1\leq i<j\leq t}{\bf Re}\left(a_i\cdot (Q_{r_j}C(y_i))(w_i),
a_j\cdot (Q_{r_j}C(y_j))(w_j)\right)\right| \geq
&& \nonumber \\
\left|c_0^2\sum_{i=1}^t|a_i|^2 - 2\|\Ea(C)\|\frac{\varepsilon_0}{4}
\right|\geq \frac{c_0^2}{2}.&&
\end{eqnarray}

Let us choose from the sequence $\{y_i\}_{i\in \N}$ a subsequence
$\{y_{l(i)}\}$ such that both closures of the subspaces $W$ and $C(W)$
spanned by the vectors $\{w_{l(i)}\}_{i\in \N}$ and
$\{Q_{r_{l(i)}}C(y_{l(i)})(w_{l(i)})\}_{i\in \N}$ respectively,
have infinite codimension. Due to the inequalities (\ref{eq50i})
and (\ref{eq50j}) the map
$$S|_W:W\to C(W), \qquad S(w_{l(i)})=Q_{r_{l(i)}}C(y_{l(i)})(w_{l(i)}),$$
is bounded and bijective. Hence, by theorem III.11 of \cite{RS}
the map $S|_W$ is an isomorphism, and we can extend it to an
invertible operator $S:H\to H$ by choosing any isomorphism between
orthogonal complements of the spaces $W$ and $C(W)$.

Since $X$ is a compact space then there exists a point $x_0\in X$
such that for any open neighbourhood $U_{x_0}$ of the point $x_0$
there are infinitely many members of the subsequence
$\{y_{l(i)}\}_{i\in \N}$ lying in $U_{x_0}$. For any open in $IM$--topology
neighbourhood $U_{\frac{1}{d}, -S, C(x_0)}$ of the operator $C(x_0)$, $d\in \N$,
we denote by $U_{x_0}^d\subset X$ the open neighbourhood
$$U_{x_0}^d:= C^{-1}\left(U_{\frac{1}{d}, -S, C(x_0)}\right)$$
of the point $x_0\in X$.

Let $z_d:= y_{l(i_d)}\in U_{x_0}^d\cap \{y_{l(i)}\}_{i\in \N}$
which exists by the choice of the point $x_0\in X$.
Then there exists a sequence of invertible operators $G_d\in GL(H)$
such that
\begin{eqnarray}\label{eq50k}
\frac{1}{d} > \left\| G_d\left(-S + Q_{r_{l(i_d)}}C(z_d)\right) -
\left(-S + Q_{r_{l(i_d)}}C(x_0)\right) \right\|\geq && \nonumber \\
\left\|\left( G_d\left(-S + Q_{r_{l(i_d)}}C(z_d)\right) -
\left(-S + Q_{r_{l(i_d)}}C(x_0)\right)\right)(w_{l(i_d)})\right\| = &&
\nonumber \\
\left\|\left(-S + Q_{r_{l(i_d)}}C(x_0)\right)(w_{l(i_d)})\right\|. &&
\end{eqnarray}
Hence, for sufficiently large $d$
\begin{equation}\label{eq50l}
\|C(x_0)(w_{l(i_d)})\|\geq \frac{1}{2}\|S(w_{l(i_d)})\|\geq \frac{c_0}{2}
\end{equation}
But the inequality (\ref{eq50l}) contradicts the condition that
$C(x_0)$ is a compact operator.
\qed
\vskip 0.2 cm

\begin{cor}
The uniform topology is stronger than the $IM$--topology.
\end{cor}
\proof Let
$X=\{0\}\cup \bigcup_{i=1}^{\infty}\{\frac{ 1}{i}\} \subset \R$
and $\Aa= C(X)$. We shall construct an $\Aa$--operator
$K\in \Ka(l_2(\Aa))$ such that the map
$\Da(K):X\to \Ka(H)^{IM}\stackrel{\id}{\longrightarrow}
\Ka(H)^s \subset \Ba(H)^s$ is not continuous in the uniform
topology but by the theorem \ref{thm8} is continuous in the $IM$--topology.

Let us define continuous functions $\varphi_i:X \to \C$, $i\in \N$,
by the following rule
$$\varphi_i\left(\frac{1}{i}\right)=1, \quad \varphi_i\left(\frac{1}{j}\right)=0, \quad for \quad j\ne i, \quad \varphi_i(0)=0.$$
We define the $\Aa$--operator $K:l_2(\Aa)\to \l_2(\Aa)$ by the following formula
\begin{equation}\label{eqAA}
K(\xi)=\left(\sum_{i=1}^{\infty}\varphi_i\xi_i, 0, 0, \ldots \right),
\end{equation}
where $\xi=(\xi_1,\xi_2,\ldots)$ is an element of the $l_2(\Aa)$,
$\xi_i\in \Aa$, $i=1\ldots, \infty$.
Then $K\in \Ka(l_2(\Aa))$ and hence the map $\Da(K):X\to
\Ka(H)^{IM}\stackrel{\id}{\longrightarrow} \Ka(H)^s \subset \Ba(H)^s$
is continuous in the $IM$--topology.

By definition (\ref{eqAA}) of the $\Aa$--operator $K$ we have
$$\Da(K)(0) = 0 \in \Ba(H).$$
But
$$\left\|\Da(K)\left(\frac{1}{i}\right)\right\| \geq
\left\|\Da(K)\left(\frac{1}{i}\right)(e_i)\right\|=1,$$
where $e_i \in H$, $=1,\ldots, \infty$ is the standard basis of $H$.
That means that the map $\Da(K):X\to \Ka(H)^{IM}\stackrel{\id}{\longrightarrow}
\Ka(H)^s \subset \Ba(H)^s$ is not continuous in the uniform topology.
\qed
\vskip 0.2 cm

Now, let us discuss the representing space for
$\K$-theory introduced in the paper \cite{AtSe-03e}.
In this paper M.~Atiyah and G.~Segal have considered locally trivial
bundles $P\to X$ whose fibers $P_x={\bf P}(H)$ are the projective space
of a separable infinite dimensional complex Hilbert space $H$ and structural
group is the projective unitary group $\Pa\Ua(H)^{c.o}$,
with the compact-open topology. With the aim to define
a twisted \K-theory they need to replace fiber $P_x$ by
a representing space for \K-theory such that the structural
group $\Pa\Ua(H)^{c.o}$ acts continuously on it by conjugation
(ibid. sect.3, p.12).
It is well known (see \cite{lect} and \cite{Jan}) that the space
$\Fa(H)^u$ of Fredholm operators in $H$ with the uniform topology
is a representing space for \K-theory. Unfortunately,
the unitary group $\Ua(H)^{c.o}$ (and $\Pa\Ua(H)^{c.o}$),
with the compact-open topology, does not act continuously
on $\Fa(H)^u$ by conjugation. To surmount this obstacle
M.~Atiyah and G.~Segal have suggested (ibid.) to use as
a representing space for  \K-theory the following set
$$Fred'(H) =\{(A,B) \in \Fa(H)\times \Fa(H)\ |\  AB-I\in \Ka(H)
\quad and \quad BA-I \in \Ka(H)\}$$
with the topology induced by the embedding
$$Fred'(H) \hookrightarrow \Ba(H)^{c.o}\times \Ba(H)^{c.o}\times
\Ka(H)^u\times \Ka(H)^u $$
$$(A,B) \to (A,B,AB-I,BA-I),$$
where $\Ba(H)^{c.o}$ is the space of bounded operators in $H$,
with the compact-open topology, and $\Ka(H)^u$ is the space of
compact operators in $H$, with the uniform topology.
Let $X$ be a compact space.
In this case, by Banach-Steinhaus theorem (see \cite{RS} Theorem III.9)
the continuous
maps $X\to \Ba(H)$ are the same for the compact-open and
for the strong operator topologies. Then any continuous map
$f: X \to Fred'(H)$ of compact space $X$ into $Fred'(H)$ can
be considered as a pair of continuous maps in the strong operator
topology
\begin{equation}\label{eq51}
A^f : X \to \Fa(H)^{s},
\end{equation}
\begin{equation}\label{eq52}
B^f : X \to \Fa(H)^{s}
\end{equation}
such that for any $x\in X$ the operators $A^f(x)B^f(x)-I$ and
$B^f(x)A^f(x)-I$ are compact and the maps
\begin{equation}\label{eq53}
A^fB^f-I : X \to \Ka(H)^u,
\end{equation}
\begin{equation}\label{eq54}
B^fA^f-I : X \to \Ka(H)^u
\end{equation}
are continuous in the uniform topology.

Now, we can relax the conditions (\ref{eq53}) and (\ref{eq54}) to
strictly extend the set of admitted maps $A^f$, in the following way.
The maps
\begin{equation}\label{eq55}
A^fB^f-I : X \to \Ka(H)^{IM},
\end{equation}
\begin{equation}\label{eq56}
B^fA^f-I : X \to \Ka(H)^{IM}
\end{equation}
are continuous in the $IM$-topology. Indeed, by the theorem \ref{thm8} we have
\begin{equation}\label{eq57}
\Ea(A^fB^f-I), \quad \Ea(B^fA^f-I) \in \Ka(l_2(\Aa)),
\end{equation}
By the theorems \ref{thm2} and \ref{thm3} the $\Aa$--operators
$\Ea(A^f)$ and $\Ea(B^f)$ are Fredholm $\Aa$--operators. Due to the
isomorphism (\ref{eq47}) the maps
\begin{equation}\label{eq58}
A^f : X \to \Fa(H)^F,
\end{equation}
\begin{equation}\label{eq59}
B^f : X \to \Fa(H)^F
\end{equation}
are continuous in the $F$-topology. Hence, the class of
continuous maps $X\to \Fa(H)^F$ is strictly wider than the class
of continuous maps $A^f: X\to \Fa(H)^s$ for which the conditions
(\ref{eq51}), (\ref{eq52}), (\ref{eq53}), (\ref{eq54}) hold.
Moreover, it was shown in the paper \cite{Irm} that the space
$\Fa(H)^F$ of Fredholm operators, with the $F$-topology, is a
representing space for \K-theory. Taking in account the theorem
\ref{thm6}, we conclude, that we can take the space $\Fa(H)^F$ as
a representing space for \K-theory in the construction of the twisted \K-theory.

Results of the sections 2,3 were obtained by A.S.Mishchenko
and results of the section 4 were obtained by A.A.Irmatov.
Research is partially supported by the grant of
RFBR No 02-01-00574, the grant of the support for Scientific
Schools No NSh-619.2003.1, and the grant of the foundation
"Russian Universities" No UR.04.03.009.


\end{document}